\documentclass[12pt]{article}
\usepackage{amssymb}

\newtheorem{theo}{Theorem}

\newtheorem{prop}{Proposition}

\newtheorem{lem}{Lemma}

\newtheorem{rem}{Remark}

\newtheorem{defi}{Definition}

\begin{document}
\author{ Aziz Ikemakhen \thanks{
Facult\'{e} des sciences et techniques, B.P. 549, Gueliz,
Marrakech, Morocco. \hspace{4cm} \ \  E-mail:
ikemakhen@fstg-marrakech.ac.ma } }
\title{ Parallel Spinors on Pseudo-Riemannian $Spin^c$ Manifolds}
\date{ March 14, 2005 }
\maketitle

\begin{abstract}
 We describe, by their holonomy groups,  all   simply connected  irreducible
 non-locally symmetric pseudo-Riemannian $Spin^c$
  manifolds which admit parallel spinors. So
  we generalize the Riemannian $Spin^c$ case (\cite{Mor}) and the pseudo-Riemannian
  $Spin$ one (\cite{BK}).
\end{abstract}
\textbf{Mathematics Subject Classifications:} 53C50, 53C27.\\
\textbf{Key words:} holonomy groups, pseudo-Riemannian $Spin^c$
manifolds, parallel spinors.
\section{Introduction}
In (\cite{Mor}),  Moroianu described all  simply connected
 Riemannian $Spin^c$ manifolds admitting parallel spinors.
Precisely, he showed the following result:
\begin{theo}
 A simply
connected $Spin^c$ Riemannian  manifolds $(M,g)$ admits a parallel
spinor if and only if it is isometric to the Riemannian product
$(M_1,g_1)\times (M_2,g_2) $ of a complete simply connected
k$\ddot{a}$hler manifold $(M_1,g_1)$ and a complete simply
connected $Spin$ manifold $(M_2,g_2) $ admitting a parallel
spinor.  The $Spin^c$ structure of $(M,g)$ is then the product of
the canonical $Spin^c$ structure of $(M_1,g_1)$ and the $Spin$
structure of $(M_2,g_2)$.
\end{theo}
In (\cite{BK}), Baum and Kath  characterized, by their holonomy
group, all simply connected  irreducible non-locally symmetric
pseudo-Riemannian $Spin$ manifolds admitting parallel spinors.
Precisely, they proved the following result:
\begin{theo}
Let $(M,g)$ be a simply connected  irreducible
 non-locally symmetric pseudo-Riemannian $Spin$ manifold of
 dimension $n=p+q$ and  signature $(p,q)$. We denote by $N$  the dimension of the
 space of parallel spinors on $M$.
 Then   $(M,g)$ admits
  a  parallel spinors if and only if the
 holonomy group $H$ of $M$ is (up to conjugacy in $O(p,q)$)
 one  in the following list :

$$
\begin{tabular}{||c|c||}
\hline \hline
  Holonomy group & N \\
  \hline
  \hline
  $SU(p',q') \subset SO(2p',2q')$ & 2 \\
  \hline
  $Sp(p',q') \subset SO(4p',4q')$ & $ p' + q' + 1$ \\
  \hline
  $ G_2 \subset SO(7)$ & 1 \\
  \hline
  $ G'_{2(2)} \subset SO(4,3)$ & 1 \\
  \hline
  $ G^{\mathbb{C}}_2 \subset SO(7)$ & 2 \\
  \hline
  $ Spin(7)\subset SO(8)$ & 1 \\
  \hline
  $ Spin(4,3) \subset SO(8, 8)$ & 1 \\
  \hline
  $Spin(7,\mathbb{C})\subset SO(8, 8)$ &  1 \\
  \hline
\end{tabular}
$$
$$(table \;\;1)$$
\end{theo}
Our aim is to generalize this result for the simply connected
irreducible non-locally symmetric pseudo-Riemannian $Spin^c$
manifolds. More  precisely, we show that:
\begin{theo} Let $(M,g)$ be a connected simply connected  irreducible
 non-locally symmetric  $Spin^c$ pseudo-Riemannian  manifold of
 dimension $n=p+q$ and  signature $(p,q)$. Then
 the following conditions are equivalent\\
 (i) $(M,g)$  admits
  a  parallel spinor,  \\
 (ii) either $(M,g)$ is a  $Spin$ manifold which admit a parallel
 spinor,
 or $(M,g)$ is a k$\ddot{a}$hler not Ricci-flat
 manifold,\\
 (iii) the
  holonomy group $H$  of $(M,g)$ is (up to conjugacy in
$O(p,q)$)
 one  in table 1 or $H = U(p',q')$, $p = 2p'$ and $q = 2q'$.\\
 For $H = U(p',q')$ the dimension of the space of parallel spinors on $M$ is 1.
\end{theo}
This theorem is a contribution to the resolution of the following
problem :\\
\textbf{(P)} What are the possible holonomy groups of simply
connected
 pseudo-Riemannian $Spin^c$ manifolds which admit parallel
 spinors?\\
 Some partial answers to this problem have been given by  Wang  for
 the Riemannian $Spin$ case (\cite{Wa}), by  Baum and  Kath  for the
 irreducible pseudo-Riemannian $Spin$ one (\cite{BK}) , by   Leistner
for the  Lorentzian $Spin$ one (\cite{Le, L2}), by   Moroianu for
the Riemannian $Spin^c$ one (Theorem 1), and by author for the
totally reducible pseudo-Riemannian $Spin$ one and the  Lorentzian
$Spin^c$ one (\cite{I1, I2} ). The problem remains open even
though big progress have been made because the classification of
the possible holonomy groups of pseudo-Riemannian manifolds  is
not yet made with the exception of the irreducible case made by
Berger (\cite{Ber1, Ber2} ) and the case of Lorentzian manifolds
made by  B\'{e}rard Bergery and the author (\cite{BI}). By De
Rham-Wu's splitting theorem, the problem \textbf{(P)} can be
reduced to the case of the indecomposable pseudo-Riemannian
manifolds (its holonomy representation  do not leave any
non-degenerate proper subspace ). But the general classification
remains extremely difficult, because some indecomposable but non
irreducible manifolds exist. I.e. its holonomy representation
leaves a degenerate proper subspace but its do not leave any
non-degenerate proper subspace. In this article we are content
only with studying the irreducible case that is a particular case
of the indecomposable one.\\
In paragraph 2 of this paper we define the group $Spin^{c}(p,q)$
and its spin representation. We also define the
$Spin^{c}$-structure on pseudo-Riemannian manifolds and its
associated spinor bundle. In paragraph 2 we give an algebraic
characterization to the pseudo-Riemannian $Spin^c$ manifolds which
admit parallel spinors and  we prove Theorem 3.
\section{ Spinor representations and  $Spin^{c}$- bundles }
\subsection{   $Spin^{c}(p,q)$ groups }
Let $<.,.>_{p,q}$ be the ordinary scalar product of signature
$(p,q)$ on $\mathbb{R}^m$ $(m=p+q)$. Let $Cl_{p,q}$ be the
Clifford algebra of $\mathbb{R}^{p,q}:=(\mathbb{R}^m
,<.,.>_{p,q})$ and $\mathbb{C}l_{p,q}$ its complexification. We
denote by $\cdot$ the Clifford multiplication of
$\mathbb{C}l_{p,q}$. $\mathbb{C}l_{p,q}$ contains the groups
$$\mathbb{S}^1:=\{z \in \mathbb{C};  \parallel z
\parallel=1  \}$$ and
$$Spin(p,q):=\{X_1 \cdot...\cdot X_{2k};  \;\;
<X_i,X_i>_{p,q} \; = \pm1;\;\; k \geq 0  \}.$$ Since $\mathbb{S}^1
\cap Spin(p,q)=\{-1,1\}$, we define the group $Spin^{c}(p,q)$ by
$$ Spin^{c}(p,q)= Spin(p,q)\cdot \mathbb{S}^1
 = Spin(p,q)\times_{\mathbb{Z}_2} \mathbb{S}^1 .$$ Consequently,
the elements of $ Spin^{c}(p,q)$ are the
 classes $[g,z]$ of pairs $(g,z) \\ \in Spin(p,q)\times
\mathbb{S}^1 $, under  the equivalence relation $(g,z)\sim
(-g,-z)$. The following suites are exact (see \cite{LM}):
$$
\displaystyle 1\rightarrow \mathbb{Z}_2 \rightarrow Spin(p,q)
\stackrel {\lambda}{\longrightarrow} SO(p,q)\rightarrow 1
$$

$$
\displaystyle 1\rightarrow \mathbb{Z}_2 \rightarrow Spin^{c}(p,q)
\stackrel {\xi}{\longrightarrow} SO(p,q)\times \mathbb{S}^1
\rightarrow 1,
$$
where $\lambda(g)(x)= g \cdot x \cdot g^{-1}$ for $ x\in
\mathbb{R}^m $  and $ \xi ([g,z])=
(\lambda(g),z^2). $\\
Let $(e_i)_{1\leq i\leq m }$ be an  orthonormal basis of
$\mathbb{R}^{p,q}$ ( $<e_i,e_j> = \varepsilon_i \delta_{ij },$
$\varepsilon_i = -1 $ for $1\leq i\leq p $  and  $\varepsilon_i =
+1 $  for  $1+p\leq i\leq m $ ). The Lie algebras of $Spin(p,q)$
and $Spin^{c}(p,q)$ are respectively
$$ spin(p,q):= \{ e_i \cdot e_j \; ; 1 \leq i < j \leq m \} $$
and
$$
 spin^{c}(p,q):= spin(p,q) \oplus \texttt{i}\mathbb{R}.
 $$
The derivative  of  $\xi$ is a  Lie algebra isomorphism and it is
given by
$$
\xi_\ast ( e_i \cdot e_j , \texttt{i}t)= (\alpha_\ast ( e_i \cdot
e_j) , \texttt{i}t) = ( 2 E_{ij} , 2\texttt{i}t ),
$$
 where $ E_{ij}= -\varepsilon_j D_{ij}+ \varepsilon_i D_{ji}$ and $
 D_{ij}$ is the standard basis of  $gl(m,\mathbb{R})$ with the
 $(i,j)$-component equal  1 and all other zero.

\subsection{$Spin^{c}$ representations }
Let $ U= \left(
\begin{array}{cc}
 0 &  i   \\
 i &  0
\end{array}
\right), \;  V= \left(
\begin{array}{cc}
 0 &  -1   \\
 1 &  0
\end{array}
\right), \; E= \left(
\begin{array}{cc}
1 &  0   \\
 0 &  1
\end{array}
\right) \; , \; T= \left(
\begin{array}{cc}
 -1 &  0   \\
 0 &  1
\end{array}
\right) $, and $\mathbb{C}(2^n)$ the complex algebra consisting of
$2^n \times 2^n$-matrices. It is well know that the Clifford
algebra $\mathbb{C}l_{p,q}$ is isomorphic to $\mathbb{C}(2^n)$ if
$m =p + q $ is even and to $\mathbb{C}(2^n) \oplus
\mathbb{C}(2^n)$ if $m$ is odd. Some natural isomorphisms  are
defined like follows (see \cite{BK}): \\
 In case $m  = 2n$ is even, we define
$ \Phi _{p,q} : \mathbb{C}l_{p,q} \rightarrow \mathbb{C}(2^n) $ by
$$ \Phi _{p,q}
(e_{2j-1})   =  \tau_{2j-1}  E \otimes ... \otimes E \otimes U
\otimes T \otimes ...\otimes T
 $$
\begin{equation}  \label{ eqn : 1}
 \;\;\;\; \;\; \Phi _{p,q} (e_{2j})   =   \tau_{2j}  E \otimes ... \otimes E \otimes V \otimes
  \underbrace{T \otimes ... \otimes T}_{\mbox{(j-1)-times}},
 \end{equation}
where $\tau_j =\mathrm{i} $ if $\varepsilon_j  = -1$ and $\tau_j
=1$ if
$\varepsilon_j  = 1.$\\
In case  $m =2n +1 $ is odd, $ \Phi _{p,q} :  \mathbb{C}l_{p,q}
\rightarrow \mathbb{C}(2^n) \oplus \mathbb{C}(2^n)$  is defined by

 $$
 \Phi _{p,q}(e_k)   = (\Phi _{p,q-1}(e_k),  \Phi
 _{p,q-1}(e_k)),\;\;
 k= 1, ..., m-1;
 $$
 \begin{equation}  \label{ eqn : 2}
\Phi_{p,q}(e_m)=(\texttt{i}T\otimes...\otimes T , -\texttt{i}
T\otimes...\otimes T ).
\end{equation}
This yields representations of the spin group $Spin (p,q)$ in case
$m$ even by restriction and in case $m$ odd by restriction and
projection onto the first component.
 the  module space of $Spin(p,q)$- representation is
  $ \Delta_{p, q}= \mathbb{C}^{2^n}$.
The   Clifford multiplication is defined by
 \begin{equation}  \label{ eqn : 3}
X \cdot u := \Phi_{p,q}(X) (u) \;\;\; \mbox{ for} \;\;\;  X \in
\mathbb{C}^m \;\;\; \mbox{and}  \;\;\;  u \in \Delta_{p, q}.
 \end{equation}
A usual basis  of $  \Delta_{p, q}$ is the following :
$u(\nu_n,...,\nu_1):=u(\nu_n) \otimes...\otimes u(\nu_1);
 \;\;   \nu_j =\pm1,$
 where \\

$u(1)=  \left(%
\begin{array}{c}
  1\\
  0
  \end{array}
  \right)   \;\; \hbox{and} \;\;
  u(-1)=  \left(%
\begin{array}{c}
  0\\
  1
  \end{array}
  \right ) \in \mathbb{C}^2.$ \\
  The spin representation of the group $Spin(p,q)$ extends to a
   $Spin^{c}(p,q)$- representation by :

\begin{equation}\label{4}
 \Phi_{p,q}([g,z])(v)= z \; \Phi_{p,q}(g) (v),
\end{equation}
for $v \in \Delta_{p,q}$ and $ [g,z] \in Spin^{c}(p,q)$. Therefore
 $\Delta_{p,q}$ becomes the  module space of
 $Spin^{c}(p,q)$- representation (see \cite{Frid}).\\
There exists a hermitian inner product $<.,.>_{\Delta}$ on the
spinor module $\Delta_{p,q}$ defined by
$$
<v,w>_{\Delta}:=  \texttt{i}^{\frac{p(p-1)}{2}}(e_1 \cdot ...\cdot
e_p\cdot v,w ) ; \;\; \hbox{for} \;\;  v, w  \in \Delta_{p,q},
$$
where $ (z,z')=\sum_{i=1}^{2^n} z_i \cdot \overline{z'_i}$ is the
standard  hermitian  product on  $\mathbb{C}^{2^n}$.\\
$<.,.>_{\Delta}$ satisfies the following properties :
\begin{equation}  \label{ eqn : 5}
<X\cdot v,w>_{\Delta} = (-1)^{p+1} < v,X\cdot w>_{\Delta},
\end{equation}
for $X \in \mathbb{C}^m$.
\subsection{Spinor bundles}

Let $(M,g)$\ be a connected pseudo- Riemannian manifold of signature $%
(p,q)$. And let  $P_{SO(p,q)}$\ denote the bundle of  oriented
positively frames on M.
\begin{defi}
A structure $Spin^{c}$ on  $(M,g)$ is the data of  a
$\mathbb{S}^1$-principal bundle  $P_{\mathbb{S}^1}$ over $M$ and a
$\xi$-reduction $(P_{Spin^c(p,q)},\Lambda)$ of the product
$(SO(p,q)\times \mathbb{S}^1)$-principal bundle  $P_{SO(p,q)}
\widetilde{\times} P_{\mathbb{S}^1}$.
i.e. $\Lambda: P_{Spin^c(p,q)} \rightarrow (P_{SO(p,q)} \widetilde{\times} P_{\mathbb{S}^1})$ is a 2-fold covering verifying:\\
i) $P_{Spin^c(p,q)}$ is a $Spin^c(p,q)$-principal bundle over $M$,\\
ii) $\forall u \in P_{Spin^c(p,q)}$,  $\forall a \in
Spin^\mathbb{c}(p,q)$,
$$ \Lambda(ua)= \Lambda(a)\xi(a). $$
A structure $Spin$ on  $(M,g)$ is the data of  a
$\lambda$-reduction $P_{Spin(p,q)}$ of $P_{SO(p,q)}$.
\end{defi}
We note  that $(M,g)$ carries a $Spin^c$-structure if and only if
the second Stiefl-Whitney class of $M$, $w_2(M) \in H^2(M,
\mathbb{Z})$ is the $ \mathbb{Z}_2$ reduction of an integral class
$ z \in H^2(M, \mathbb{Z}_2 )$ (\cite{LM}, \cite{Frid} ).\\
\textbf{Example 1} \\
 Every pseudo-Riemannian  $Spin$ manifold is canonically a
$Spin^c$ manifold. The $Spin^c$- manifold is obtained as
$$
P_{Spin^c(p,q)} = P_{Spin(p,q)}\times_{\mathbb{Z}_2}\mathbb{S}^1,
$$
where $P_{Spin(p,q)}$ is the Spin-bundle and $\mathbb{Z}_2$ acts
diagonally by $(-1,-1)$.
\\
\textbf{Example 2} \\  Any  irreducible  pseudo-Riemannian not
Ricci-flat k$\ddot{a}$hler manifold is
canonically a $Spin^c$ manifold.\\
Indeed The holonomy group $H$ of $(M,g)$ is $U(p',q')$, where
$(p,q)= (2p',2q')$ is the signature of $(M,g)$ . Then
 $P_{SO(p,q)}$ is reduced to the holonomy $U(p',q')$-principal bundle
 $P_{U(p',q')}$. Moreover there exists an
 $<.,.>_{p,q}$-orthogonal almost complex structure  $J$
 which we can imbed $U(p',q')$ in
 $SO(p,q)$ by
 $$
\begin{array}{ccc}
  i : U(p',q') & \hookrightarrow & SO(p,q) \\
  A + \texttt{i} B =((a_{kl})_{1\leq k ,l \leq m} + \texttt{i}(b_{kl})_{1\leq k ,l \leq m}) &
   \rightarrow & \left(%
   \begin{array}{c} \left(%
\begin{array}{cc}
  a_{kl} & b_{kl} \\
 -b_{kl} & a_{kl} \\
\end{array}%
\right)
\end{array}%
\right)_{1\leq k ,l \leq m} \\
\end{array},
$$
and $( e_k  , J e_k )_{k=1,...,p'+q'}$ is an orthogonal basis of
$\mathbb{R}^{p,q}$. \\
 We consider the homomorphism
$$
\begin{array}{cccc}
  \alpha: & U(p',q') & \hookrightarrow & SO(p,q)\times \mathbb{S}^1 \\
   & C & \rightarrow & (i(C),\det(C)) \\
\end{array}
$$
The proper values of every element $C \in U(p',q')$ is in
$\mathbb{S}^1$ and
$$
 \cos 2 \theta + \varepsilon_k \sin 2\theta \;\; e_k \cdot J e_k
= \varepsilon_k (\cos \theta \;\; e_k  +   \sin \theta \;\; J e_k
)\cdot (\sin \theta \;\; e_k  -   \cos \theta \; \; J e_k ),
$$
where \,\, $ \varepsilon_k = \; <e_k,e_k>_{p,q} $.  Then the
following homomorphism is well defined :
$$
\begin{array}{cccc}
 \widetilde{\alpha}: & U(p',q') & \hookrightarrow & Spin^c(p,q) \\
   & C & \rightarrow & \displaystyle \prod_{k=1}^m (
    \cos 2 \theta_k + \varepsilon_k \sin 2\theta_k \; \; e_k \cdot Je_k) \times
   \displaystyle e^{ \frac{\texttt{i}}{2} \sum \theta_k},
\end{array}
$$
where  $\displaystyle e^{\texttt{i} \theta_k}$, $k=1,...,m$, are
the proper values of $C$. And it is easy to verifies that the
following diagram commutes
\begin{equation}\label{6}
    \begin{array}{lll}
&  & Spin^{c}c(p,q) \\
& \stackrel{\stackrel{\sim }{\alpha }}{\nearrow } & \downarrow \xi  \\
U(p^{\prime },q^{\prime }) & \stackrel{\alpha }{\longrightarrow }
& SO(p,q).
\end{array}
\end{equation}
Consequently,
$$
P_{Spin^c(p,q)} =
P_{U(p',q')}\times_{\widetilde{\alpha}}Spin^c(p,q).
$$
Now, let  denote by $S := P_{Spin(p,q)}
\times_{\Phi_{p,q}}\Delta_{p,q}$ the spinor bundle associated to
the $Spin^c$-structure $P_{Spin(p,q)}$. The  Clifford
multiplication given by (\ref{ eqn : 3}) defines a Clifford
multiplication on $S$:
$$
\begin{array}{ccc}
  TM \otimes S = (P_{Spin(p,q)} \times_{\Phi_{p,q}}\mathbb{R}^m )\otimes (P_{Spin(p,q)}
\times_{\Phi_{p,q}}\Delta^{\pm}_{p,q} )& \rightarrow & S \\
 ( X\otimes \psi ) = [q,x]\otimes [q,v]& \rightarrow & [q,x\cdot v]=:X\cdot\psi. \\
\end{array}
$$
Since the   scalar product $<.,.>_{\Delta}$ is $Spin_0
^c(p,q)$-invariant, it defines  a scalar product on $S$ by :
$$ <\psi,\psi_1>_{\Delta}= <v,v_1>_{\Delta},\;\; \hbox{for} \;\; \psi = [q,v]
\;\; \hbox{and} \;\; \psi_1 = [q,v_1] \in \Gamma(S).
$$
According to (\ref{ eqn : 5}), it is then easy to verify that
\begin{equation}  \label{ eqn : 6}
 <X\cdot \psi,\psi_1>_{\Delta} = (-1)^{p+1} <\psi,X\cdot
 \psi_1>_{\Delta},
\end{equation}
for $ X \in \Gamma(M)$ and   $\psi, \psi_1 \in \Gamma(S)$.\\
 Now,
as in the  Riemannian case ( see \cite{Frid}),
 if $(M, g)$ is a $Spin^{c}$ pseudo- Riemannian manifold, every connection form $A:
TP_{\mathbb{S}^1} \rightarrow \mathrm{i} \mathbb{R}$  on the
$\mathbb{S}^1$- bundle $P_{\mathbb{S}^1}$ defines \\( together
with the Levi-Civita $D$ of $(M, g)$ ) a covariant derivative
$\nabla^A$ on the spinor bundle $S$, called
the spinor derivative  associated to  $(M,g, S, P_{\mathbb{S}^1} , A)$.\\
Henceforth, a $Spin^{c}$ pseudo- Riemannian manifold will be the
data of a set \\ $(M,g, S, P_{\mathbb{S}^1}, A)$, where $(M,g)$ is
an oriented connected pseudo- Riemannian manifold, $S$ is a
$Spin^{c}$ structure, $P_{\mathbb{S}^1}$ is the
$\mathbb{S}^1$-principal  bundle  over $M$ and $A$ is a connection
form on $P_{\mathbb{S}^1}$. Using (7) and by the same proof in the
Riemannian case ( see \cite{Frid}), we conclude that
\begin{prop} $ \forall \; X, Y \in \Gamma(M) $  and
$\forall \; \psi, \psi_1 \in \Gamma(S)$,
\begin{equation}  \label{ eqn : 7}
  \nabla^A _Y  (X\cdot
\psi) =X\cdot \nabla^A_Y  (\psi) + D_Y X \cdot \psi.
\end{equation}
\begin{equation}  \label{ eqn : 8}
 X< \psi,\psi_1>_{\Delta} =  <\nabla^A_X \psi,
 \psi_1>_{\Delta}+ <\psi,\nabla^A_X \psi_1>_{\Delta}.
\end{equation}
\end{prop}
Let us denoted by $F_A := \mathrm{i} w $ the  curvature form of
$A$, seen as an imaginary-valued  2-form on $M$,  by $R$ and $Ric$
respectively  the curvature
  and the Ricci tensors of $(M,g)$ and by $R^A$ the curvature tensor of $\nabla^A$.
  Like  Riemannian case ( see
  \cite{Frid}), if we put $A(X):=X\lrcorner \omega $ we have
\begin{prop} For $q=(e_1,...,e_m)$ a local  section  of \ \ $P_{Spin(p,q)}$,\\
$ \forall \; X, Y \in \Gamma(M) $ \hbox{and} $\forall \; \psi \in
\Gamma(S)$,
\begin{equation}  \label{ eqn : 9}
R^A(X,Y)\psi = \frac{1}{2}\sum _{1\leq i<j\leq m}\varepsilon_i
\varepsilon_j \; g(R(X,Y) e_i , e_j ) e_i\cdot e_j\cdot \psi +
\mathrm{i} \frac{1}{2} \; \omega(X,Y)\cdot \psi,
\end{equation}
and
\begin{equation}  \label{ eqn : 10}
\sum _{1\leq i\leq m}\varepsilon_i \; R^A(X,e_i)\psi = -
\frac{1}{2} Ric(X)\cdot \psi + \mathrm{i} \frac{1}{2} A(X)
\cdot\psi.
\end{equation}
\end{prop}
\begin{rem}  According to Example 1, if $(M,g)$ is
 $Spin$
 then it is $Spin^{c}$. Moreover,
  the  auxiliary bundle $P_{\mathbb{S}^1}$ is trivial and then there exists a
global  section  $ \sigma : M \rightarrow P_{\mathbb{S}^1}$. We
choose the connection  defined by $A$ to be flat, and we denote $
\nabla^A $ by $\nabla .$  Conversely, if the auxiliary bundle
$P_{\mathbb{S}^1}$ of a $Spin^{c}$-structure is trivial, it is
canonically identified with a
 $Spin$-structure. Moreover, if the  connection $A$ is flat,
by this identification, $\nabla^A$ corresponds to the
 covariant derivative on the spinor bundle.
\end{rem}

\section{parallel Spinors}
\subsection{algebraic characterization }
It is well know that there exists a bijection between  the space
$\mathcal{PS}$ of all parallel spinors on $(M,g)$   and the space
$$ V_{\widetilde{H}} = \{v \in
\Delta_{p,q}
 ;  \;\; \Phi_{p,q}(\widetilde{H})(v):= \widetilde{H}\cdot v = v \}
 $$
of all fixed spinors of $\Delta_{p,q}$ with respect to the
holonomy group  $\widetilde{H}$ of the  connection $\nabla^A$
(\cite{BK}). If $(M,g)$ is supposed  simply connected, then
$\mathcal{PS}$ is in bijection with $$ V_{\mathcal{H}}= \{v \in
\Delta_{p,q}
 ; \;\; \widetilde{\mathcal{H}}\cdot v = 0 \},$$
where $\widetilde{\mathcal{H}}$ is  the Lie algebra of
$\widetilde{H}$. Moreover, provided with the connection led by the
Levi-Civita connection $D$ and the form connection $A$ the
holonomy group of $P \times P_{\mathbb{S}^1}$ is  $\xi
(\widetilde{H})\subset H \times H_A$, where $H$ is the holonomy
group of $(M,g)$ and $H_A$ the one of $A$ (see propositions 6.1
and 6.3 \cite{KN1}) . $H_A = \{1\}$ if $A$ is flat and $H_A =
\mathbb{S}^1$ otherwise.  With the notations introduced in
subsection 2.1, for $(B , it) \in \xi_*
(\widetilde{\mathcal{H}})$, where $\widetilde{\mathcal{H}}$ is the
Lie algebra of $\widetilde{H}$  we have
$$
\xi^{-1}_{\ast} (B,  \mathrm{i} t) = ( \lambda^{-1}_{\ast} (B),
\frac{1}{2}\mathrm{i}t ).
$$
Now if we differentiate the relation (4) at [1,1],  we get :
$$ \phi_{p,q}(C, \mathrm{i} t)(v) = \mathrm{i} t v + \phi_{p,q}(C)(v) ,$$
for $(C , it) \in  spin^c(p,q)$ and $v \in \Delta_{p,q}$.
 Then
$$ \phi_{p,q}( \xi_\ast^{-1}(B, it))(v) =
 \frac{1}{2} \mathrm{i} t v + \phi_{p,q}(\lambda^{-1}_{\ast} (B))(
 v).
$$
 We conclude that
\begin{prop}
 $(M,g)$ admits a non trivial parallel spinor  if and only
if there exists  \\ $0 \neq v \in \Delta_{p,q}$ such that

\begin{equation}\label{12}
   B \cdot v := \phi_{p,q}(\lambda^{-1}_{\ast} (B))(
 v) = -\frac{1}{2} \mathrm{i} t v,
   \;  \forall \;\; (B, it) \in \xi_\ast(\widetilde{\mathcal{H}})\subset
    \mathcal{H}\oplus \mathcal{H}_A,
\end{equation}
where  $\widetilde{\mathcal{H}}$, $\mathcal{H}$ and $\mathcal{H}_A
$ are respectively the  Lie algebras of $\widetilde{H}$, $H$ and
$H_A$.
\end{prop}

\subsection{  Proof  of  Theorem 3 }

Let $(M,g, S, P_{\mathbb{S}^1}, A)$ be a $Spin^{c}$ structure
where $(M , g)$ is  a connected simply connected  irreducible
 non-locally symmetric pseudo-Riemannian  manifold of
 dimension $n=p+q$ and  signature $(p,q)$,
 which admits a non trivial parallel spinor
$\psi$. \\ We consider the two  distributions $T$ and $E$ defined
by
$$  T_x := \{ X \in T_xM;  \;\; X\cdot \psi = 0 \},
$$
$$
 E_x = \{ X \in T_xM; \; \exists \; Y \in T_xM; \;\; X\cdot \psi = \texttt{i}
 Y\cdot\psi \},
 $$
for $x \in M$.  Since $\psi$ is parallel, By (8), $T$ and $E$ are
parallel. Since  $T$ is isotropic and the manifold $(M,g)$  is
supposed irreducible, by the holonomy principe, we have
\begin{equation}\label{13}
   T=0.
\end{equation}
Now denote by $F$ the image of the Ricci tensor:
$$ F_x := \{ Ric(X);  \;\; X \in T_xM \}.$$
Since $\psi$ is parallel, (11) shows that
\begin{equation}\label{14}
   Ric(X)\cdot \psi = \texttt{i} A(X)\cdot \psi.
\end{equation}
Then $F \subset E$. Consequently, from (13), we obtain
$$  E^\perp  \subset F^\perp = \{ Y \in TM;\; Ric(Y)=0
 \} = \{ Y \in TM;\; A(Y )  = 0 \}.$$
 $(M,g)$  is supposed irreducible, by the holonomy
principe, $E= 0$ or $E= TM$. \\
If $E= 0$, then $F=0$. This gives $Ric=0$ and $A=0$. According to
Remark 1, $(M,g)$ is Spin and $\psi$ is a parallel spinor on
$M$.\\
If $E= TM$, we have a $(1,1)$-tensor $J$ definite by
\begin{equation}\label{15}
   X \cdot \psi  = \texttt{i} J(X)\cdot
 \psi,\; \hbox{where} \; X \in TM .
\end{equation}
\begin{lem} For $X, Y \in T_xM $, if
 \ \ $( X +  \texttt{i} Y )  \cdot \psi = 0$ then \\ $ g(X,Y)
=0$ and $g(X,X) = g(Y,Y)$.
\end{lem}
\textbf{Proof of Lemma 1}. If we denote by $g^{c}$ the complex
form of $g$   then
$$ \begin{array}{lll}
 ( X +  \texttt{i} Y )  \cdot ( X +  \texttt{i} Y )  \cdot \psi & =
 & g^{c}( X +  \texttt{i} Y , X +  \texttt{i} Y  )  \psi \\
   & = & ( g(X,X) - g(Y,Y) + 2 i g(X,Y) ) \psi  = 0. \\
\end{array}
$$
Since $\psi$ is non trivial we obtain the lemma.\\
 Lemma 1 implies that $J$ defines an
orthogonal  almost complex structure on $M$. Moreover, from (8)
and (15) we obtain $J$ is parallel, since $\psi$ is parallel. In
consequence, $(M,g)$ is a
k$\ddot{a}$hler  manifold.\\
Now if $(M,g)$ is a k$\ddot{a}$hler  manifold then there exists a
canonical $Spin^c$ structure of $(M,g)$. And
from Remark 1, the following conditions are equivalent \\
(a) $(M,g)$ is not $Spin$, \\
(b)  $H_A = \mathbb{S}^1 $, \\
(c)  $(M,g)$  is not Ricci-flat, \\
(d) $H = U(p',q')$.\\
Then  the equivalence between (i) and (ii) of Theorem 3  are
proved. And from Theorem 2, we have the  equivalence between (ii)
and (iii). To finish the proof of Theorem 3, it remains to show
for $H = U(p',q')$ that the dimension of parallel spinors on $M$
is $1$. For this, we remark that if we reduce the principle
bundles  $P_{Spin^c(p,q)}$ and $P_{SO(p,q)} \widetilde{\times}
P_{\mathbb{S}^1}$ to their holonomy bundles the diagram  (6)
becomes
\begin{equation}\label{}
    \begin{array}{lll}
&  & \widetilde{H} \\
& \stackrel{\stackrel{\sim }{\alpha }}{\nearrow } & \downarrow \xi  \\
 U(p^{\prime },q^{\prime }) & \stackrel{\alpha
}{\longrightarrow } & \alpha ( U(p^{\prime },q^{\prime })),
\end{array}
\end{equation}
and the holonomy group of $P_{SO(p,q)} \widetilde{\times}
P_{\mathbb{S}^1}$ is exactly $\alpha ( U(p^{\prime },q^{\prime }))
= \xi (\widetilde{H} )$.
 However
$U(p',q')=SU(p',q') \times U_{\mathbb{S}^1}$, where
$$
U_{\mathbb{S}^1} = \{ \;\; \left(%
\begin{tabular}{llll}
$\lambda $ & 0 & $\cdots $ & 0 \\
0 & 1 & $\ddots $ & $\vdots $ \\
$\vdots $ & $\ddots $ & $\ddots $ & 0 \\
0 & $\cdots $ & 0 & 1
\end{tabular}
\right); \;\;  \lambda \in \mathbb{S}^1 \;\;
 \},
$$
$u(p',q') = su(p',q') \oplus u_{\mathbb{S}^1}$  where
$u_{\mathbb{S}^1} \simeq \texttt{i}\mathbb{R} $ \ \ is the Lie
algebra of $U_{\mathbb{S}^1}$. If we consider the imbedding
$$
\begin{array}{ccc}
  i : u(p',q') & \hookrightarrow & so(2p',2q') \\
  A+ \texttt{i}B =((a_{kl})_{1\leq k ,l \leq n} + \texttt{i}(b_{kl})_{1\leq k ,l \leq n}) &
   \rightarrow & \left(%
   \begin{array}{c} \left(%
\begin{array}{cc}
  a_{kl} & b_{kl} \\
 -b_{kl} & a_{kl} \\
\end{array}%
\right)
\end{array}%
\right)_{1\leq k ,l \leq n} \\
\end{array}
$$
$\alpha_*(u_{\mathbb{S}^1})$ is generated  by $(E_{12}, \mathrm{i}
)$. Then
$$
\begin{array}{ccc}
  \xi_* (\widetilde{\mathcal{H}}) =  \alpha_*(u(p',q')) & = & \alpha_*(su(p',q')) \oplus \alpha_*(u_{\mathbb{S}^1}) \\
  & = & i(su(p',q')) \oplus (E_{12}, \mathrm{i} )  \mathbb{R}  \\
\end{array}
$$
 From \cite{BK}, $u^+ := u(1,...,1)$ and $u^- :=
u(-1,...,-1)$ generate the space \\ $ V_{su(p',q')}= \{ v \in
\Delta_{p,q} \;
 ; \;\; su(p',q') \cdot v=0 \}
 $.  Moreover, by (1)
$$
E_{12} \cdot u^+ = \frac{1}{2} e_1 \cdot e_2 \cdot u^+ =
\frac{1}{2} \texttt{i} u^+ \; \; \hbox{and}\; \; E_{12} \cdot u^-
= -\frac{1}{2} \texttt{i} u^- .
$$
Then $u^-$ belongs to  the space
$$
 V_{u(p',q')}= \{ v \in \Delta_{p,q} \;
 ; \;\; B \cdot v + \frac{1}{2} \mathrm{i} t v =0,
   \;  \forall \;\; (B, it) \in \xi_\ast(\widetilde{\mathcal{H}})
   \},
 $$
 and it is easy to verify that $u^-$ generates $V_{u(p',q')}$.
This completes  the proof of Theorem 3.
\begin{rem} By Theorem 3, we deduce the results of Moroianau made in the
riemannian case (see \cite{Mor}).
\end{rem}


\begin{thebibliography}{9}


\bibitem{BK}   Baum,  H. and  Kath, I.:  \textit{Parallel spinors and holonomy
groups on pseudo- Riemannian spin manifolds}. Ann. Glog. Anal.
Geom., 17:1-17 (1999).
\bibitem{BI}   B\'{e}rard Bergery, L. and  Ikemakhen, A. : \textit{On the Holonomy
of Lorentzian Manifolds,} Proceeding of Sympo. in Pure Math., V
\textbf{54} (1993), Part2.
\bibitem{Ber1}  Berger, M. : \textit{Sur les groupes
d'holonomie des vari\'{e}t\'{e}s \`{a} connection affine et des
vari\'{e}t\'{e}s riemanniennes,} Bull. Soc. Math. France 83
(1955), 279-330.
\bibitem{Ber2} Berger, M. : \textit{Les espaces symétriques non
compacts,} Ann. Sci. Ecole Norm. Sup. 74 (1957), 85-177.
\bibitem{Frid}  Friedrich, Th.: \textit{Dirac Operators in
Riemannian Geometry.} Graduate Studies in Mathematics. Volume 25.
AMS Providence, Rhode Island 2000.
\bibitem{I1}   Ikemakhen, A.: \textit{Groupes d'holonomie
et spineurs paralleles sur les vari\'{e}t\'{e}s
pseudo-Riemanniennes compl\`{e}tement r\'{e}ductibles}, C. R.
Acad.
 Sci. Paris, Ser. I 339 (2004) 203-208.
 \bibitem{I2}   Ikemakhen, A.: \textit{Parallel Spinors
  on Lorentzian $Spin^c$ Manifolds}, Pr\'{e}print Universit\'{e}
  Cadi-Ayyad (2005).
\bibitem{LM}   Lawson, H.B. and  Michelsohn, M.-L.: \textit{Spin geometry} Princeton
Univ.Press 1989.
\bibitem{Le}  Leistner, Th.: \textit{Holonomy and Parallel
Spinors in  Lorentzian Geometry}, PhD thesis, Humbold-University
of Berlin, (2003).
\bibitem{L2} Leistner, Th.:  \textit{Towards a classification of Lorentzian holonomy
groups. Part II: Semisimple, non-simple weak-Berger algeabras},
arXiv:math.DG$\diagup$0309274 v1, (2003).
\bibitem{Mor} Moroianu, A.: \textit{Parallel and Killing Spinors on
$Spin^c$ Manifolds}. Commun. Math. Phys. 187, 417-427 (1997).
\bibitem{KN1}   Kobayashi, S. and  Nomizu, K. : \textit{Foundations of
differentiable geometry}, Vol.1, Interscience, Wiley, New York
(1963).
\bibitem{Wa}  Wang, M.Y.: \textit{Parallel spinors and parallel forms}, Ann. Global Anal.
Geom. 7 (1989) 1, 59-68.


\end{thebibliography}
\end{document}